\documentclass[12pt,reqno,intlim]{amsart}

\usepackage{amssymb,amsmath}

\usepackage[arrow,matrix,curve]{xy}

\newdir{ >}{{}*!/-5pt/\dir{>}}

\newcommand{\nc}{\newcommand}

\nc{\bC}{\bold{C}} \nc{\bN}{\Bbb{N}} \nc{\cF}{\mathcal{F}}
\nc{\cE}{\mathcal{E}} \nc{\cR}{\mathcal{R}} \nc{\cM}{\mathcal{M}}
\nc{\al}{\alpha} \nc{\bt}{\beta} \nc{\gm}{\gamma} \nc{\dl}{\delta}
\nc{\om}{\omega} \nc{\sg}{\sigma} \nc{\Sg}{\Sigma} \nc{\vf}{\varphi}
\nc{\ve}{\varepsilon} \nc{\os}{\overset} \nc{\ol}{\overline}
\nc{\ul}{\underline} \nc{\us}{\underset} \nc{\sbs}{\subset}
\nc{\bsl}{\backslash} \nc{\Ra}{\Rightarrow}
\nc{\lra}{\longrightarrow} \nc{\all}{\allowdisplaybreaks}

\nc{\Codes}{\operatorname{{\bold{Codes}}}}
\nc{\RegMono}{\operatorname{\mathcal{R}{\rm{eg}\mathcal{M}{\rm{ono}\!}}}}
\nc{\RegEpi}{\operatorname{\mathcal{R}{\rm{eg}\mathcal{E}{\rm{pi}\!}}}}
\nc{\Mn}{\operatorname{\mathcal{M}{\rm{ono}\!}}}
\nc{\Ep}{\operatorname{\mathcal{E}{\rm{pi}\!}}}
\nc{\Rg}{\operatorname{\mathcal{R}{\rm{eg}\!}}}
\nc{\Ob}{\operatorname{Ob\!}}

\numberwithin{equation}{section}

\theoremstyle{definition}

\theoremstyle{remark}

\begin{document}

\title[]
{A remark on normal closures in free products of groups}

\author{Dali Zangurashvili}

\maketitle

\begin{abstract}
It is shown that a nontrivial normal subgroup $N$ of a group $G$ is a free factor of the $N$'s normal closure in the $G$'s free product with arbitrary nontrivial groups.

\bigskip

\vskip+1mm
\noindent{\bf 2020  Mathematics Subject Classification}: 20E06, 20E07; 20E15; 20J15. 
\end{abstract}

\section{Introduction}

The notion of the free product of groups is one of the classical and frequently used notions in group theory. Let $I$ be a set, and $G_i$ be a nontrivial group, for any $i\in I$.  Recall that the free product $\amalg_{i\in I}G_i$ of groups $G_i$ is defined as a group $G$ containing all the groups $G_i$ and such that each element $g$ of $G$ can be uniquely represented as a reduced word, i.e., as the expression $1$ or an expression
\begin{equation}
a_1a_2...a_k
\end{equation} 
\noindent with $k>0$, $a_i\in G_{\beta(i)}$, $a_i\neq 1$ ($1\leq i\leq k$) and $\beta(i)\neq \beta(i+1)$ ($1\leq i\leq k-1$). The free product $\amalg_{i\in I}G_i$ exists and is unique up to isomorphism (see, e.g., \cite{K1}). As is well-known, it is the coproduct of the groups $G_i$ in the category of groups (see, e.g., \cite{M1}).

It is obvious that the groups $G_i$ are not normal in their free product. In view of this fact, the problem of studying their normal closures (i.e., the smallest normal subgroups $H_i$ of $\amalg_{i\in I}G_i$ containing $G_i$'s) naturally arises. It is easy to notice that, for any $i\in I$, the intersection of $H_i$ with the subgroup generated by the subgroups $G_j$ with $j\neq i$ is trivial (this subgroup is clearly the free product $\amalg_{j\neq i, j\in I} G_j$). This implies that the quotient group $G/H_i$ is isomorphic to $\amalg_{j\neq i, j\in I} G_j$ (see, e.g., \cite{LAL}).

In the present note, we deal with arbitrary normal subgroups of groups $G_i$, and show that these subgroups are free factors of their normal closures in the free product $\amalg_{i\in I}G_i$.

\section{Normal closures in free products of groups}

Recall that, by the Kurosh subgroup theorem \cite{K}, for any subgroup $H$ of the free product $G$ of groups $G_i$ ($i\in I$), there are a set $J$ and subgroups $F$ and $A_j$ ($j\in J$) of $G$ such that $F$ is a free group, each $A_j$ is conjugate to a subgroup of some $G_i$, and one has the representation 
\begin{equation}
H=F\amalg (\amalg_{j\in J} A_j).
\end{equation}
Later, MacLane proved that there is representation (2.1) with $J=I$ and 
\begin{equation}
A_i=\amalg_{D} B_{i, D},
\end{equation}
\noindent for each $i\in I$, where $D$ ranges over the $(H,G_i)$-double cosets in $G$, and
\begin{equation}
B_{i, D}=H\cap sG_{i}s^{-1},
\end{equation}
\noindent for some element $s$ in the double coset $D$ \cite{M}. To this end, MacLane introduced the notions of an $i$-system of representatives (for $i\in I$ and the subgroup $H$ of the free product $\amalg_{i\in I}G_i$), and of a uniform Schreier system of representatives (for the subgroup $H$). An $i$- system of representatives is a function $R_i$ which assigns to each left coset $C$ of $H$ a representative $R_i(C)$ of $C$, with $R_i(H)=1$, such that $R_i(Ca)\in R_i(C)G_i$, for any $a\in G_i$. A uniform Schreier system of representatives is defined as
a collection $(R_i)_{i\in I}$ of $i$-systems that satisfies a certain analog of the Schreier condition from the familiar algebraic proof of the fact that
every subgroup of a free group is free (for details, we refer the reader to \cite{M}; see also \cite{H}, \cite{K1}).
MacLane proved the existence of a uniform Schreier system of representatives, for any $H$.

Moreover, MacLane showed that, for an $i$-system $R_i$ of representatives and a $(H,G_i)$-double coset $D$, there is a unique element $s\in D$ such that $sG_i$ contains all $R_i(C)$ for $C\subseteq D$, and such that either $s=1$ or $\omega(s)\neq i$; here $\omega(g)$ denotes the index $\beta(k)$ of the group to which the last "syllable" in the reduced form (1.1) of $g$ belongs. These $s$'s were called $i$-double representatives. 

MacLane proved that $i$-double representatives $s$'s render the representation (2.1) (equipped with equalities (2.2) and (2.3)) valid, provided that the collection of $i$-systems $R_i$ of representatives is a uniform Schreier system. 

As Kurosh noted in \cite{K1} (p. 458), in that case and for such $s$'s, all nontrivial intersections $H\cap G_i$ participate in the decomposition 
\begin{equation}
H=F\amalg (\amalg_{i\in I} (\amalg_{D} B_{i, D})).
\end{equation}
Indeed, consider a $(H,G_i)$-double coset $D$ containing $H$. Any left coset $C$ contained in $D$ has the form $Ha$ with $a\in G_i$. Then, by the definition of an $i$-system of representatives, $R_i(C)=R_i(Ha)\in R_i(H)G_i=G_i$. By the uniqueness of $s$, $s=1$. This implies that $B_{i.D}=H\cap G_i$.
\vskip+2mm
 
The above-mentioned facts imply the following
\vskip+2mm
\textbf{Theorem.}
Let $G_i$ be nontrivial groups $(i\in I)$. Let $N$ be a nontrivial normal subgroup of $G_i$, for some $i$, and $H$ be its normal closure in the free product $G=\amalg_{j\in I} G_j$. Then there is a subgroup $K$ of $H$ such that $H=N\amalg K$; the subgroup $K$ is unique up to isomorphism.

\begin{proof}
The embedding $\mu_i:G_i\rightarrowtail \amalg_{j\in I} G_j$ is obviously a split monomorphism. Let $\varepsilon_i$ be its left inverse. Since the subgroup $H$ is generated by the elements of the form $gag^{-1}$ with $g\in G$ and $a\in N$, one has $\varepsilon_i(H)\subseteq N$. This implies that $H\cap G_i=N$ (this equality immediately follows also from Theorem 3.1 of \cite{Z} since any split monomorphism is a codescent morphism \cite{JT1}). Therefore, $N$ participates in decomposition (2.4), for relevant $s$'s. 

The second claim of the statement follows from the observation given in Section 1.


\end{proof}

The author gratefully acknowledges the financial support from Shota Rustaveli National Science Foundation of Georgia (Ref.: STEM-22-1601).
\vskip+5mm

\vskip+2mm

\textit{Author's address:}

\textit{Andrea Razmadze Mathematical Institute of Tbilisi State University}, 

\textit{6 Aleksidze Str., Lane II, Tbilisi 0193, Georgia;}

\textit{e-mail: dali.zangurashvili@tsu.ge}

\end{document}